# On DeTurck uniqueness theorems for Ricci tensor

Sergey Stepanov

*Finance University, Leningradsky Prospect, 49-55, 125468 Moscow, Russian Federation*

___

**Abstract**

In Riemannian geometry the *prescribed Ricci curvature problem* is as follows: given a smooth manifold $M$ and a symmetric 2-tensor $r$, construct a metric on $M$ whose Ricci tensor equals $r$. In particular, DeTurck and Koiso proved the following celebrated result: the Ricci curvature uniquely determines the Levi-Civita connection on any compact Einstein manifold with non-negative section curvature.

In the present paper we generalize the result of DeTurck and Koiso for a Riemannian manifold with non-negative section curvature. In addition, we extended our result to complete non-compact Riemannian manifolds with nonnegative sectional curvature and with finite total scalar curvature.

*MSC*: 53C20

*Keywords:* uniqueness theorem for Ricci tensor; compact and complete Riemannian manifolds; vanishing theorem

___

## 1. Introduction

The point of the paper [1] and the monograph [2; pp. 140-153] is that in certain circumstances the metric (or at last the connection) is uniquely determined by the Ricci tensor. In particular, in [1, Corollary 3.3] and [2, Theorem 5.42] anyone can read the following: Let $(M, \bar{g})$ be a compact *Einstein manifold* with non-negative section curvature and with the Ricci tensor $Ric(\bar{g}) = \bar{g}$, then another Riemannian metric $g$ on $M$ with $Ric(g) = \bar{g}$ has the same Levi-Civita connection as $\bar{g}$. We remark that this proposition is a corollary of the Eells and Sampson *vanishing theorem* for harmonic maps (see [3, p. 124]).

In the present paper we consider a compact *Riemannian manifold* $(M, \bar{g})$ with non-negative sectional curvature and with $Ric(\bar{g}) \leq \bar{g}$. Under these conditions, we prove that if $g$ is another Riemannian metric on $M$ with the Ricci tensor $Ric(g) = \bar{g}$, then $g$ and $\bar{g}$ have the same Levi-Civita connection. Furthermore, if the full holonomy group Hol($\bar{g}$) is irreducible then the metric $g = C\bar{g}$ for some constant $C > 0$. In turn, it is well known that $Ric(\bar{g}) = Ric(C\bar{g})$. This proposition is a corollary of our vanishing theorem for harmonic maps which was announced in [4]. We extend the above scheme to show that if $(M, \bar{g})$ is a non-compact manifold $(M, \bar{g})$ with non-negative sectional curvature and with the Ricci tensor $Ric(\bar{g}) \leq \bar{g}$ then there is no complete Riemannian metric $g$ such that its Ricci tensor $Ric(g) = \bar{g}$ and its total scalar curvature $s_g(M)$ is finite. This proposition is a corollary of the Schoen and Yau *vanishing theorem* for harmonic maps of complete non-compact Riemannian manifolds (see [7]).

Our propositions complement the result of the paper [1] and the monograph [2].

## 2. Harmonic maps

For the discussion of harmonic maps we will follow Eells and Sampson [3]. Let $(M, g)$ and $(\overline{M}, \bar{g})$ be two Riemannian manifolds with the Levi-Civita connections $\nabla := \nabla(g)$ and $\overline{\nabla} := \nabla(\bar{g})$, and $f : (M, g) \to (\overline{M}, \bar{g})$ be a smooth map. The *energy density* of $f$ is defined as the scalar function

$$e(f) = 2^{-1} \| df \|^2 \qquad (1.1)$$

where $\| df \|^2$ is the squared norm of the differential of $f$ with respect to metric on the bundle $T^*M \otimes f^*T\overline{M}$. Then the *total energy* of $f$ is obtained by integrating the energy density $e(f)$ over $M$

$$E(f) = \int_M e(f) dv_g \qquad (1.2)$$

where $dv_g$ denotes the measure on $(M, g)$ induced by the metric $g$. If $f$ is of class $C^2$ and $E(f) < +\infty$, and $f$ is an extremum of the *Dirichlet energy functional* $E(f)$, then $f$ is called a *harmonic map* and satisfies the *Euler-Lagrange equation*

$$\text{trace}_g D df = 0 \qquad (1.3)$$

where $D$ is the connection in the bundle $T^*M \otimes f^*T\overline{M}$ induced from the Levi-Civita connections $\nabla$ and $\overline{\nabla}$ of $(M, g)$ and $(\overline{M}, \bar{g})$, respectively.

For any harmonic $f:(M,g)\to(\overline{M},\overline{g})$ we have the *Weitzenböck formula* (see [3])

$$\Delta e(f) = Q(f) + \|D\,df\|^2 \quad (1.4)$$

where $\Delta$ is the Laplace–Beltrami operator $\Delta = div\,\nabla$ and

$$Q(f) = g(Ric, f^*\overline{g}) - trace_g(trace_g(f^*\overline{Riem})) \quad (1.5)$$

where $Ric = Ric(g)$ is the Ricci tensor of $(M, g)$ and $\overline{Riem}$ is the Riemannian curvature tensor of $(\overline{M},\overline{g})$. Let the inequality $\overline{sec} \leq 0$ be satisfied anywhere on $(\overline{M},\overline{g})$ and the inequality $Ric \geq 0$ be satisfied anywhere on compact $(M, g)$, then $Q(f)$ is non-negative everywhere on $M$. Since our hypothesis implies that the left hand side of (1.3) is non-negative, then using the *Hopf's lemma* (see [5, pp. 30-31]), one can verify that $e(f)$ is constant. In this case, from (1.4) we obtain $D\,df = 0$. In this case, $f$ is totally *geodesic map* (see [6]). Now we can formulate the following *vanishing theorem* on harmonic maps. Namely, if $f:(M,g)\to(\overline{M},\overline{g})$ is any harmonic mapping between a compact Riemannian manifold $(M, g)$ with the Ricci tensor $Ric \geq 0$ and a Riemannian manifold $(\overline{M},\overline{g})$ with the sectional curvature $\overline{sec} \leq 0$ then $f$ is *totally geodesic* and has constant *energy density* $e(f)$. Furthermore, if there is at least one point of $M$ at which its Ricci curvature $Ric > 0$, then every harmonic map $f:(M,g)\to(\overline{M},\overline{g})$ is constant (see [3, p. 124]).

In turn, Schoen and Yau have showed in [7] that $\sqrt{e(f)}$ is subharmonic function on $(M, g)$ if $Q(f) \geq 0$. On other hand, Yau has proved in his other paper [8] that every non-negative $L^2$-integrable subharmonic function on a complete Riemannian manifold must be constant. Applying this to $\sqrt{e(f)}$, we conclude that $\sqrt{e(f)}$ is a constant if the total energy $E(f) < +\infty$ (see also [7]). On the other hand, every complete non-compact Riemannian manifold with nonnegative Ricci curvature has infinite volume (see [8]). In our case, we have $Ric \geq 0$ then the volume of $(M, g)$ is infinite. This forces the constant $e(f)$ to be zero and $f$ to be a constant map (see also [7]). Now we can formulate another celebrated *vanishing theorem* on harmonic maps: If the sectional curvature of $(\overline{M},\overline{g})$ is non-positive and $(M, g)$ is a complete non-compact manifold with $Ric \geq 0$, then any harmonic map $f : (M, g) \to (\overline{M}, \overline{g})$ with the finite energy $E(f)$ is a constant map (see [7]; [9, p. 116]). We remark that in the original paper [7] the manifold $(\overline{M}, \overline{g})$ was assumed to be compact. However, this assumption is superfluous (see [9, p. 116]).

## 3. The main theorem

If we consider the manifold $M$ with two Riemannian metrics $g$ and $\overline{g}$ then the identity mapping $\mathrm{Id} : (M, g) \to (M, \overline{g})$ is harmonic if and only if the deformation tensor $T = \overline{\nabla} - \nabla$ is a section of the tensor bundle $TM \otimes S_0^2 M$, because in this case the Euler-Lagrange equation (1.3) has the form $\mathrm{trace}_g T = 0$ (see [1; 3]). In particular, if $(M, g)$ is a manifold of strictly positive Ricci $Ric$ curvature, then $\mathrm{Id} : (M, g) \to (M, Ric)$ is a harmonic map (see [1]). Next we can formulate and prove the following

**Theorem 1**. *Let $(M, \overline{g})$ be a compact Riemannian manifold with the sectional curvature $\overline{\mathrm{sec}} \geq 0$ and with the Ricci tensor $\overline{Ric} \leq \overline{g}$. If $g$ is another Riemannian metric on $M$ with the Ricci tensor $Ric = \overline{g}$, then $g$ and $\overline{g}$ have the same Levi-Civita connection. Furthermore, if the full holonomy group $\mathrm{Hol}(\overline{g})$ of $(M, \overline{g})$ is irreducible then $Ric = \overline{Ric}$.*

**Proof.** With the above assumptions, we have $Ric = \overline{g} > 0$, then the identity map $\mathrm{Id} : (M, g) \to (M, \overline{g})$ is *harmonic*. In this case, we have $e(f) = \frac{1}{2} s$ for the energy density $e(f)$ of the harmonic identity map $\mathrm{Id} : (M, g) \to (M, \overline{g})$ and the scalar curvature $s = \mathrm{trace}_g Ric$ of the Riemannian manifold $(M, g)$ (see [1]; [2; p. 152]). Therefore, $s$ satisfies the Weitzenböck formula (1.4) which has the following form (see [1]):

$$\tfrac{1}{2} \Delta s = Q(f) + \| D \mathrm{d} f \|^2 \qquad (3.1)$$

where $Q(f) = g^{ik} g^{jl} (\overline{g}_{ij} \overline{g}_{kl} - \overline{R}_{ijkl})$ and $\| D \mathrm{d} f \|^2 = g^{ij} g^{kl} \overline{g}_{pq} T_{ik}^p T_{jl}^q \geq 0$ for local components $g_{ij}, \overline{g}_{kl}, \overline{R}_{ijkl}$ and $T_{kl}^i$ of metric tensors $g$ and $\overline{g}$, the Riemannian

curvature tensor $\overline{Riem}$ and the deformation tensor $T$, respectively. On the other hand, we have the identity (see [2, p. 436] and [10])

$$\left(\overline{g}_{ij}\overline{R}_{kl} - \overline{R}_{ijkl}\right)\varphi^{ik}\varphi^{jl} = \sum_{i<j}\overline{sec}\left(\overline{e}_i,\overline{e}_j\right)\left(\overline{\lambda}_i - \overline{\lambda}_j\right)^2 \quad (3.2)$$

where $\varphi$ is any smooth symmetric tensor field such that $\varphi(\overline{e}_i,\overline{e}_j) = \overline{\lambda}_i\delta_{ij}$ for the Kronecker delta $\delta_{ij}$ and some orthonormal basis $\{\overline{e}_1,...,\overline{e}_n\}$ at any point $x \in M$. Then equation (3.1) can be rewritten in the form

$$\tfrac{1}{2}\Delta s = \sum_{i<j}\overline{sec}\left(\overline{e}_i,\overline{e}_j\right)\left(\overline{\lambda}_i - \overline{\lambda}_j\right)^2 + g^{ik}g^{jl}\overline{g}_{ij}\left(\overline{g}_{kl} - \overline{R}_{kl}\right) + \|T\|^2 \quad (3.3)$$

where $g(\overline{e}_i,\overline{e}_j) = \overline{\lambda}_i\delta_{ij}$. We remark that under the stated assumptions the right side of (3.3) is non-negative, since then $\Delta s \geq 0$. Therefore, the scalar curvature $s$ is a *positive subharmonic function* on $(M, g)$. If $(M,\overline{g})$ is a compact Riemannian manifold, then using the *Hopf's lemma* (see [5, pp. 30-31]), one can verify that $s = \text{const}$. In this case, from (3.3) we obtain $T = 0$. Then $g$ and $\overline{g}$ have the same Levi-Civita connection, i.e. $\overline{\nabla}g = 0$. Furthermore, if the full holonomy group Hol($\overline{g}$) of $(M,\overline{g})$ is irreducible then the metric $g = C\overline{g}$ for some constant $C > 0$ (see [2, pp. 282; 285-287]). In this case, we have the identity $Ric = \overline{Ric}$ because $Ric(\overline{g}) = Ric(C\overline{g})$ for some positive constant $C$ (see [2, pp. 44; 152]). QED.

## 4. Two vanishing theorems

In [1] the following non-existence theorem was proved: Let $(M,\overline{g})$ be a compact Riemannian manifold with all sectional curvature less then $(\overline{n}-1)^{-1}$. Then there is no Riemannian metric $g$ on $M$ such that its Ricci tensor $Ric = \overline{g}$. In turn, we can formulate and prove this proposition in the following form.

**Theorem 2**. *Let $(M, \overline{g})$ be a compact Riemannian manifold with nonnegative section curvatures and with the Ricci tensor $\overline{Ric} \leq \overline{g}$. If in addition there is at least one point of $M$ at which the Ricci tensor $\overline{Ric} < \overline{g}$, then there is no Riemannian metric $g$ on $M$ such that its Ricci tensor $Ric = \overline{g}$.*

**Proof.** Let $M$ be a compact manifold. We may assume that $M$ is oriented by taking the twofold covering of $M$ if necessary. Then by *Green's theorem* (see [5, pp. 31-33]) we obtain from (3.1) the following identity

$$\int_M Q(f)dv_g + \int_M \|T\|^2 dv_g = 0. \qquad (3.4)$$

If the inequalities $\overline{sec} \geq 0$ and $\overline{Ric} \leq \bar{g}$ are satisfied and there is a one point $x$ of $M$ in which $\overline{Ric} < \bar{g}$ then the inequality $\int_M Q(f)dv_g > 0$ holds. This inequality contradicts the equation (3.4). In this case, the harmonic mapping $f$ must be constant. QED.

In the case when $M$ is a non-compact manifold, we can prove the following

**Theorem 3.** *Let $(M, \bar{g})$ be a non-compact Riemannian manifold with the section curvature $\overline{sec} \geq 0$ and with the Ricci tensor $\overline{Ric} \leq \bar{g}$. Then there is no complete Riemannian metric $g$ on $(M, \bar{g})$ such that its Ricci tensor $Ric = \bar{g}$ and its total scalar curvature $s(M)$ is finite.*

**Proof.** Let $(M, \bar{g})$ be a non-compact Riemannian manifold with the section curvature $\overline{sec} \geq 0$ and with the Ricci tensor $\overline{Ric} \leq \bar{g}$, then $Q(f)$ is non-negative everywhere on $M$. If we assume that there is complete Riemannian metric $g$ on $(M, \bar{g})$ such that its Ricci tensor $Ric = \bar{g} > 0$, then the volume of $(M, g)$ is infinite (see [8]). Moreover, we have $e(f) = \frac{1}{2} s$ for the energy density $e(f)$ of the harmonic identity map $\mathrm{Id}: (M, g) \to (M, \bar{g})$ and the scalar curvature $s = \mathrm{trace}_g \, Ric$ of the Riemannian manifold $(M, g)$. In this case, $\sqrt{s}$ is a strictly positive subharmonic function on a complete Riemannian manifold $(M, g)$ of infinite volume (see [7]). In addition, if we suppose that the total scalar curvature $s(M) = \int_M s\,dv_g < +\infty$, then $s_g$ must be zero (see [7]; [11, p. 262]). On the other hand, according to the condition of our theorem the scalar curvature $s = \mathrm{trace}_g \bar{g} >$

0 and hence there is no complete Riemannian metric $g$ on non-compact $(M, \bar{g})$ such that its Ricci tensor $Ric = \bar{g}$. QED.